\documentclass{article}
\usepackage{xcolor,caption,subcaption,graphicx}
\usepackage{amsmath}

\newcommand\eq[1] {(\ref{#1})}

\newcommand\sect[1] {\ref{sec:#1}}
\newcommand\figref[1] {\ref{#1}}

\newcommand\labsect[1] {\label{sec:#1}}

\newcommand{\bfm}[1]{\mbox{\boldmath ${#1}$}}
\newcommand{\nonum}{\nonumber \\}
\newcommand{\beqa}{\begin{eqnarray}}
\newcommand{\eeqa}[1]{\label{#1}\end{eqnarray}}
\newcommand{\beq}{\begin{equation}}
\newcommand{\eeq}[1]{\label{#1}\end{equation}}
\newcommand{\bbm}{\begin{bmatrix}}
\newcommand{\ebm}{\end{bmatrix}}
\newcommand{\bpm}{\begin{pmatrix}}
\newcommand{\epm}{\end{pmatrix}}

\newcommand{\Real}{\mathop{\rm Re}\nolimits}

\newcommand{\lang}{\langle}
\newcommand{\rang}{\rangle}

\newcommand{\ov}[1]{\overline{#1}}

\newcommand{\wt}[1]{\widetilde{#1}}
\newcommand{\wh}[1]{\widehat{#1}}

\newcommand{\Ga}{\alpha}
\newcommand{\Gb}{\beta}
\newcommand{\Gd}{\delta}

\newcommand{\Gve}{\varepsilon}

\newcommand{\Gc}{\chi}

\newcommand{\Gs}{\sigma}

\newcommand{\BGc}{\bfm\chi}

\newcommand{\BGs}{\bfm\sigma}

\newcommand{\BGG}{\bfm\Gamma}

\newcommand{\CE}{{\cal E}}

\newcommand{\CH}{{\cal H}}

\newcommand{\CJ}{{\cal J}}

\newcommand{\CP}{{\cal P}}

\newcommand{\CU}{{\cal U}}

\def\Be{{\bf e}}

\def\Bk{{\bf k}}

\def\Bp{{\bf p}}

\def\Bu{{\bf u}}

\def\Bx{{\bf x}}

\def\BA{{\bf A}}

\def\BC{{\bf C}}

\def\BE{{\bf E}}

\def\BI{{\bf I}}
\def\BJ{{\bf J}}

\def\BL{{\bf L}}

\def\BP{{\bf P}}
\def\BQ{{\bf Q}}
\def\BR{{\bf R}}
\def\BS{{\bf S}}
\def\BT{{\bf T}}

\def\BX{{\bf X}}

\def \RR {{\mathbb R}}

\def \ba {\begin{array}}
\def \ea {\end{array}}

\def \refe #1.{(\ref{#1})}
\def \reff #1.{figure~\ref{#1}}
\def \refs #1.{section~\ref{#1}}
\def \refss #1.{subsection~\ref{#1}}
\def \refD #1.{Definition~\ref{#1}}
\def \refT #1.{Theorem~\ref{#1}}
\def \refL #1.{Lemma~\ref{#1}}
\def \refC #1.{Corollary~\ref{#1}}
\def \refP #1.{Proposition~\ref{#1}}
\def \refR #1.{Remark~\ref{#1}}
\def \refE #1.{Example~\ref{#1}}
\def \refN #1.{Notation~\ref{#1}}

\usepackage{graphics}
\usepackage{amssymb}

\begin{document}
\vspace{-1in}
\title{Substitution of subspace collections with nonorthogonal
subspaces to accelerate Fast Fourier Transform methods applied to
conducting composites}%
\label{chap:substitution}
\index{subspace collections!substitution}
\author{Graeme W. Milton}
\date{\small{Department of Mathematics, University of Utah, Salt Lake City, UT 84112, USA}}
\maketitle
\begin{abstract}We show the power of the algebra of subspace collections
developed in Chapter 7 of \cite{Milton:2016:ETC}. Specifically we accelerate the Fast Fourier Transform schemes of
Moulinec and Suquet (\cite{Moulinec:1994:FNM}, \cite{Moulinec:1998:NMC})
and \cite{Eyre:1999:FNS}, for computing the fields and effective tensor in a conducting periodic medium
by substituting a subspace collection with nonorthogonal
subspaces inside one with orthogonal subspaces. This can be done when the effective conductivity as a function of the
conductivity $\Gs_1$ of the inclusion phase (with the matrix phase  conductivity set to $1$) has its singularities
confined to an interval $[-\Gb,-\Ga]$ of the negative real $\Gs_1$ axis. Numerical results of Moulinec and Suquet
show accelerated convergence for the model example of a square array of squares at $25\%$ volume fraction.
For other problems we show how $Q^*_C$-convex functions%
\index{QCon@$Q^*_C$-convex function}
can be used to restrict the region where singularities
of the effective tensor as a function of the component tensors might be found.
\end{abstract}
\section{Introduction}
\setcounter{equation}{0}

Subspace collections with orthogonal subspaces%
\index{subspace collections!orthogonal}
obviously have a lot of
relevance to physical problems. Here we give an example to show that
subspace collections with nonorthogonal subspaces%
\index{subspace collections!nonorthogonal}
are an important tool
in analysis, reaching beyond their connection with rational functions of
several complex variables. This chapter is mostly self-contained. While it uses nonorthogonal subspace collections, it is not necessary for the reader
to have read Chapter 7 of \cite{Milton:2016:ETC}. However, it is recommended for the reader to look at Chapters 1 and 2 in that book.

\begin{figure}[t]
\centering
\includegraphics[width=0.75\textwidth]{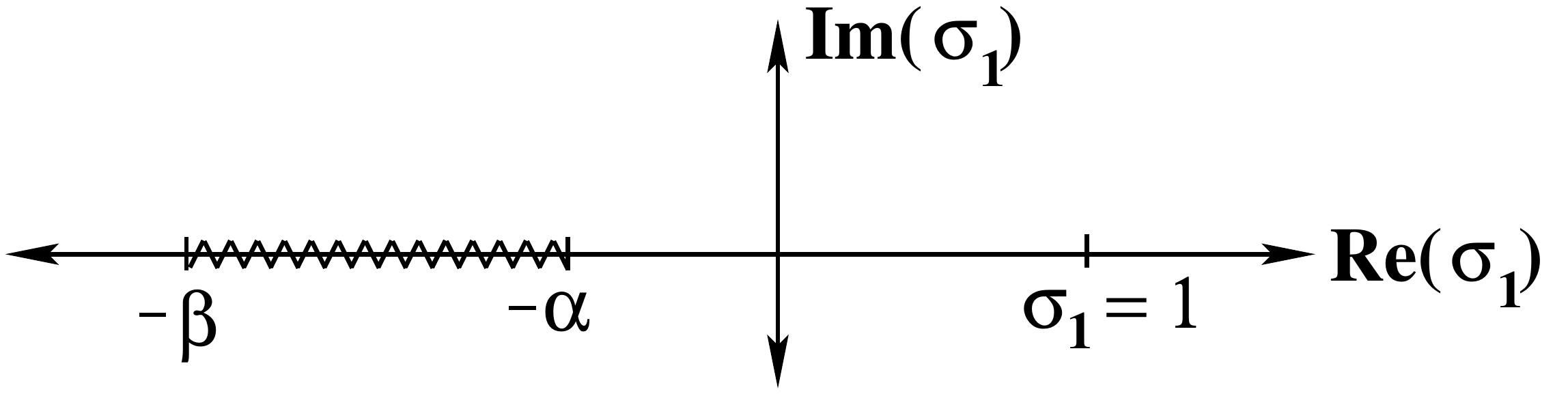}
\caption[The function $\Gs_*(\Gs_1)$ has all its singularities on the negative real $\Gs_1$-axis, between $\Gs_1=-\Ga$ and $\Gs_1=-\Gb$.]{We assume that in the complex $\Gs_1$-plane, the function $\Gs_*(\Gs_1)$ has all its singularities
on the negative real $\Gs_1$-axis, between $\Gs_1=-\Ga$ and $\Gs_1=-\Gb$. Here the zig-zag line denotes a possible branch cut.}\label{Jfig1}
\end{figure}

The specific problem we consider is a two-component conducting medium with
conductivities $\Gs_1$ and $\Gs_2$, phase $1$ being the inclusion phase.
The dimension is $d=2$ or $d=3$. Without loss of generality we can take
$\Gs_2 = 1$. Assume for simplicity that the composite has square $(d=2)$
or cubic $(d=3)$ symmetry so that the effective conductivity tensor is
isotropic, taking the form $\Gs_*\BI$. \cite{Bergman:1978:DCC} made the pioneering observation
that the function $\Gs_*(\Gs_1)$ is an analytic function%
\index{effective!conductivity analyticity}
of $\Gs_1$ and
has all its singularities on the negative real axis. He made some assumptions which
were not valid \cite{Milton:1979:TST}. These assumptions could be circumvented
by approximating the composite by a large resistor network (\cite{Milton:1981:BCP}).
Subsequently Golden and Papanicolaou (\cite{Golden:1983:BEP},
\cite{Golden:1985:BEP}) gave a rigorous proof of the analytic properties. Here we make
the additional assumption (only true for some geometries) that the function $\Gs_*(\Gs_1)$ has all its
singularities confined to an interval $[-\Gb,-\Ga]$ on the negative real
axis, $\Ga>0$ and $\Gb>\Ga$ (see \figref{Jfig1}). This information results in tighter bounds
on the conductivity and complex conductivity, or equivalently the complex dielectric constant,
(\cite{Bruno:1991:ECS}; \cite{Sawicz:1995:BCP}; \cite{Golden:1998:IMS}), and these
bounds may be inverted to yield information about $\Ga$ and $\Gb$ from experimental measurements (\cite{Orum:2012:RIS}).
By contrast, here our goal is to utilize the information
to improve the speed of convergence of Fast Fourier Transform methods for computing
the fields in composites and their associated effective tensors. These methods
were first introduced by Moulinec and Suquet%
\index{Fast Fourier Transform methods!Moulinec--Suquet scheme}
(\cite{Moulinec:1994:FNM}, \cite{Moulinec:1998:NMC}):
see \cite{Moulinec:2014:CTA} for a recent review. One important application has been to viscoplastic polycrystals
(\cite{Lebensohn:2001:NSM}, \cite{Lee:2011:MVM}).
For materials with a linear response, previous significant advances in the acceleration of these schemes were made by  \cite{Eyre:1999:FNS}%
\index{Fast Fourier Transform methods!Eyre--Milton scheme}
(see also the generalization in
section 14.9 of \cite{Milton:2002:TOC}, and in particular equation (14.38)),
and by \cite{Willot:2014:FBS}%
\index{Fast Fourier Transform methods!Willot scheme}
and \cite{Willot:2015:FBS}.
The singularities in the interval $[-\Gb,-\Ga]$ could be poles or could
be a branch cut%
\index{branch cut}
(if the inclusion has sharp corners, or if there is some
randomness in the geometry).  Note it is not
only the effective conductivity which has this analytic form but also
the electric field $\BE(\Bx,\Gs_1)$ and current field $\BJ(\Bx,\Gs_1)$
at fixed $\Bx$, and fixed applied field, with $\Gs_1$ varying. If we
didn't know anything about $\Ga$ and $\Gb$ there could potentially be
singularities anywhere along the negative real $\Gs_1$ axis. Let's first look
at the case $\Ga=0$, $\Gb=\infty$.

The original Fast Fourier Transform scheme of Moulinec and Suquet (\cite{Moulinec:1994:FNM}, \cite{Moulinec:1998:NMC})
\index{Fast Fourier Transform methods!Moulinec--Suquet scheme}
is based on series expansions, such as (2.35) in  \cite{Milton:2016:ETC}, which for conductivity takes the form
\beq \BGs_* = \Gs_0\BI+ \sum_{j=0}^\infty\BGG_0[\BGs(\Bx)-\Gs_0\BI][\BGG_1(\BI-\BGs/\Gs_0)]^j\BGG_0,\quad \Be = \Be_0+ \sum_{j=0}^\infty[\BGG_1(\BI-\BGs/\Gs_0)]^j\Be_0, \eeq{J0-0}
for the effective conductivity $\BGs_*$, and electric field $\Be(\Bx)$, with $\Be_0$ being the applied (average) electric field. [To obtain \eq{J0-0} from (2.35) in
\cite{Milton:2016:ETC} we use the
fact that $\BGG_0(\Gs_0\BI)\BGG_1=0$.]  Their key and beautiful idea was that since the action of $\BGG_1$ is readily evaluated in Fourier space,
while the action of $(\BI-\BGs/\Gs_0)$ is readily evaluated in real space, both $\BGs_*$ and $\Be(\Bx)$ can be readily calculated from these series by
going back and forth between real and Fourier space, using Fast Fourier Transforms to do so.

\begin{figure}[t]
\centering
\includegraphics[width=0.75\textwidth]{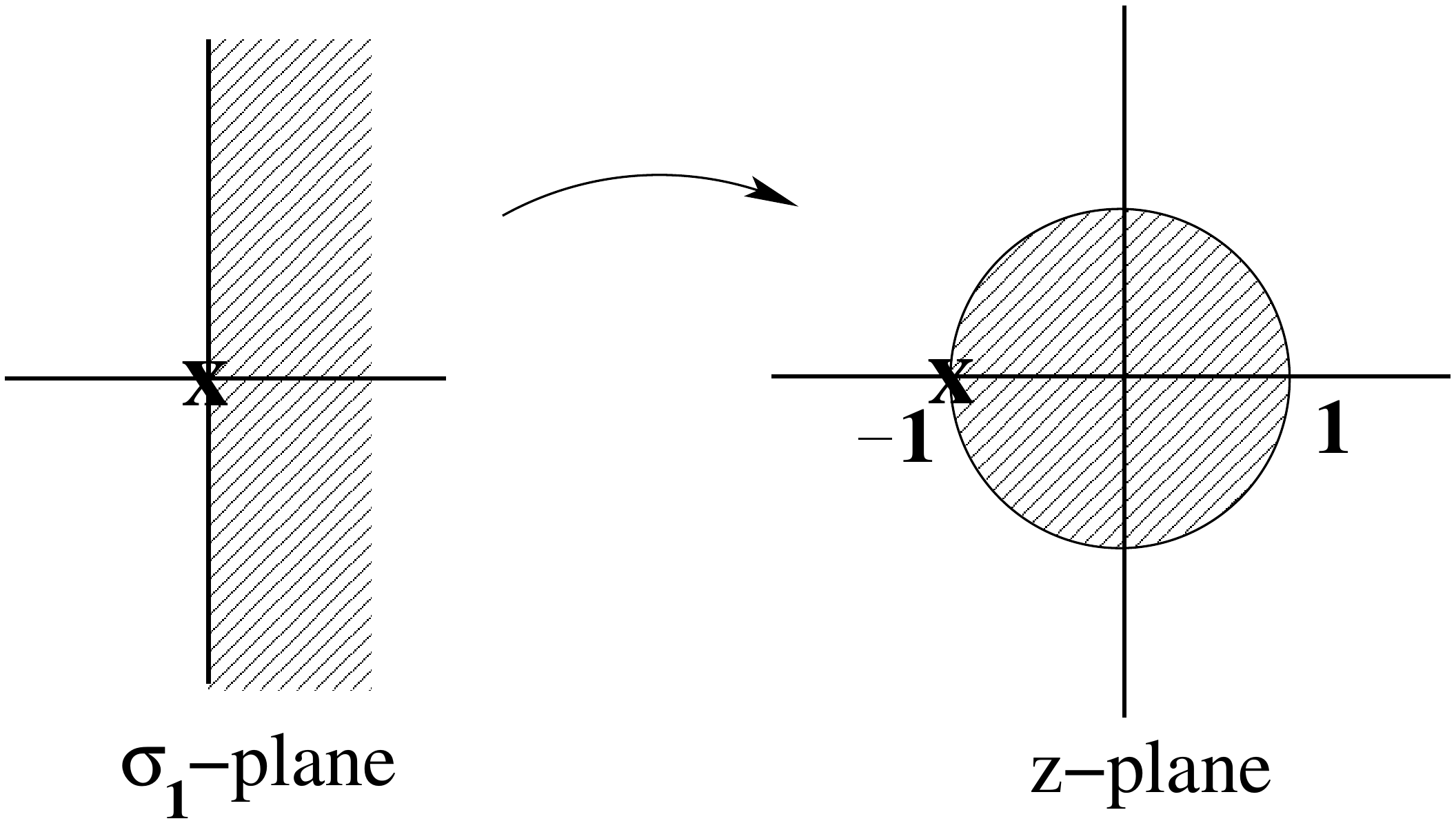}
\caption[The scheme of Moulinec and Suquet corresponds to series generated by mapping the right half of the $\Gs_1$-plane to the unit disk in the $z$-plane.]
  {If we want a series expansion which converges in the entire right half of the $\Gs_1$-plane, when a singularity is located at the point $\BX$, then we make a
fractional linear transformation%
\index{fractional linear transformation}
which takes the right half of the $\Gs_1$-plane to the unit disk in the $z$-plane, and find an series expansion%
\index{series expansions}
in powers of $z$. The
scheme of  Moulinec and Suquet (\protect\cite{Moulinec:1994:FNM}, \protect\cite{Moulinec:1998:NMC}) provides such an expansion.}\label{Jfig2}
\end{figure}

In a two-phase medium with isotropic conductivities $\Gs_1$ and $\Gs_2=1$, and with the choice of
$\Gs_0 =\frac{1}{2}(\Gs_1+\Gs_2) = \frac{1}{2}(\Gs_1+1)$ used by Moulinec and Suquet, $(\BI-\BGs/\Gs_0)$ takes the value $(1-\Gs_1)/(\Gs_1+1)$ in phase 1
and the value $-(1-\Gs_1)/(\Gs_1+1)$ in phase 2. So it is clear that this scheme gives an expansion of the form
\beq
\Gs_*/\Gs_0 = 1+\sum_{n=1}^\infty a_n\left(\frac{\Gs_1-1}{\Gs_1+1}\right)^n.
\eeq{J1}
The transformation $z = (\Gs_1-1)/(\Gs_1+1)$ maps the right half of
the complex plane to the unit disk (see \figref{Jfig2}).
So with $\Ga = 0$ the series will converge if $|z|<1$, i.e., in the
entire right half of the complex $\Gs_1-$plane, $\Real(\Gs_1)>0$. (Recall that it is the distance from the origin to the nearest singularity in the $z-$plane which
determines the radius of convergence in the $z-$plane, and hence the region of convergence in the $\Gs_1-$plane.)

\begin{figure}[!ht]
\centering
\includegraphics[width=0.75\textwidth]{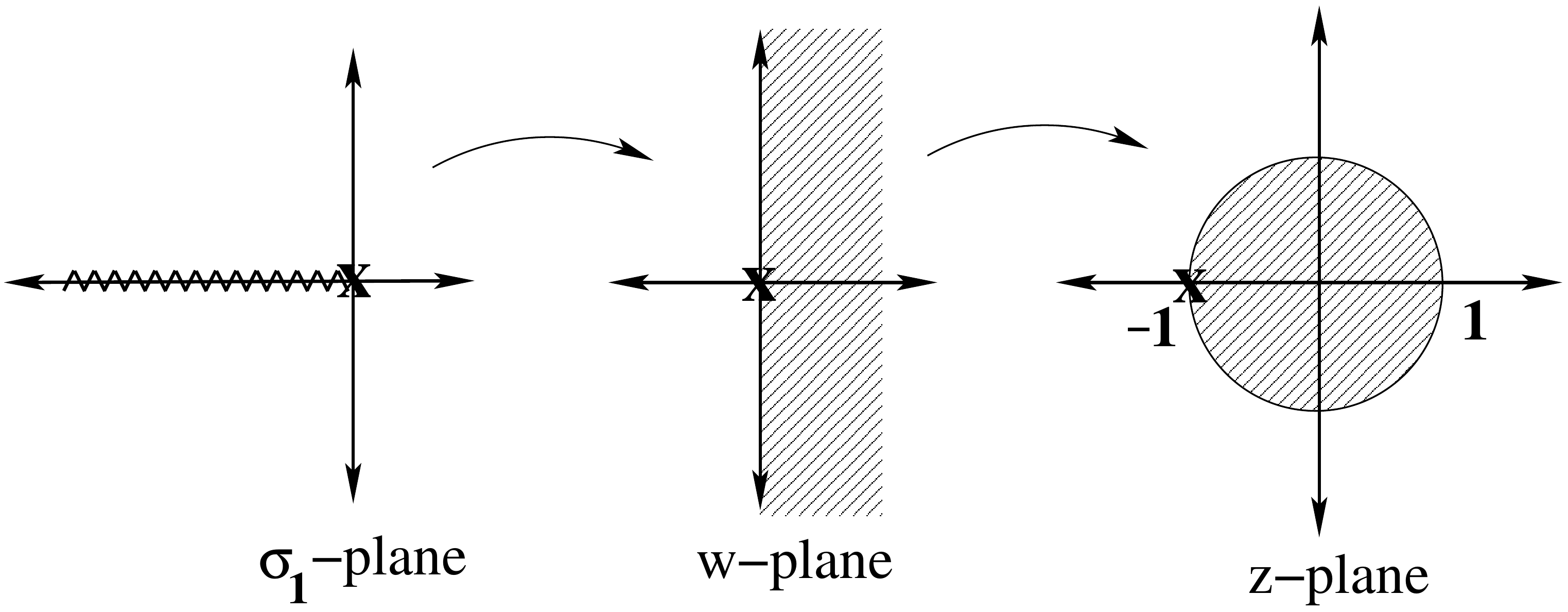}
\caption[The scheme of Eyre and Milton corresponds to series generated by mapping the cut $\Gs_1$-plane, with a cut along the negative real axis, to the unit disk in the $z$-plane.]{If we want a series expansion which converges in the entire the $\Gs_1$-plane minus the negative real $\Gs_1$-axis, then we first make a
square root transformation%
\index{square root transformation}
which maps the cut complex $\Gs_1$-plane to the right half of the $w$-plane, followed by a fractional linear transformation%
\index{fractional linear transformation}
which takes it to the unit disk in the $z$-plane, and find an expansion in powers of $z$. The
scheme of \protect\cite{Eyre:1999:FNS} provides such an expansion.}\label{Jfig3}
\end{figure}

With the accelerated scheme of \cite{Eyre:1999:FNS}%
\index{Fast Fourier Transform methods!Eyre--Milton scheme}
(see also the generalization in
section 14.9 of \cite{Milton:2002:TOC}, and in particular equation (14.38)), one has the expansion
\beq
\Gs_*/\sqrt{\Gs_1} = 1+\sum_{n=1}^\infty b_n \left(\frac{\sqrt{\Gs_1}-1}{\sqrt{\Gs_1}+1}\right)^n.
\eeq{J2}
The transformation $z = (w-1)/(w+1)$ where $w=\sqrt{\Gs_1}$ maps the
$\Gs_1$ complex plane minus the slit along the negative real axis to the
unit disk $|z|=1$ (see \figref{Jfig3}).
So with $\Ga = 0$ the series will converge if $|z|<1$, i.e., in the
entire complex plane minus the slit along the negative real axis. These arguments show that when $\Ga = 0$ the accelerated scheme
should have a larger region of convergence, and by the same line of reasoning a faster rate of convergence when both schemes converge. (Note, however, that if
$\Ga> 0$ and $\Gb<\infty$ the scheme of Moulinec and Suquet could outperform that of Eyre and Milton for small and very large values of $\Gs_1$, since the scheme of Moulinec and Suquet should then
converge for sufficiently small or sufficiently large negative real values of $\Gs_1$.)

\begin{figure}[!ht]
\centering
\includegraphics[width=0.75\textwidth]{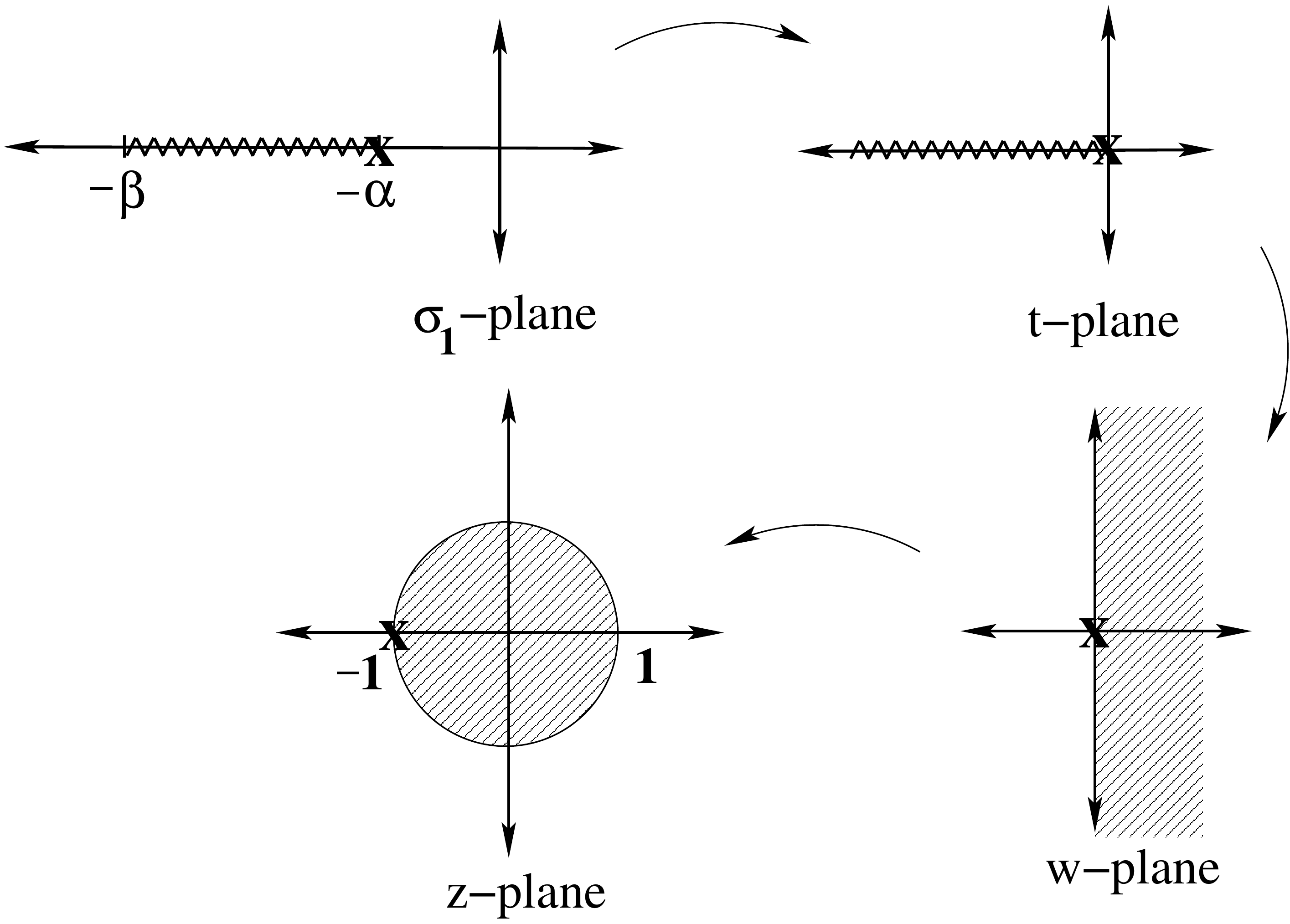}
\caption[The desired scheme should correspond to series generated by mapping the cut $\Gs_1$-plane, with a cut along the negative real $\Gs_1$-axis between $-\Gb$ and $-\Ga$, to the unit disk in the $z$-plane.]{If we want a series expansion which converges in the entire the $\Gs_1$-plane minus the cut on $\Gs_1$-axis between $-\Gb$ and $-\Ga$ then we
first use a fractional linear transformation%
\index{fractional linear transformation}
which maps this cut to the entire negative real axis in the $t$-plane, followed by a
square root transformation%
\index{square root transformation}
mapping the cut complex $t$-plane to the right half of the $w$-plane, followed by a fractional linear transformation
which takes it to the unit disk in the $z$-plane, and find an expansion in powers of $z$. The
scheme developed in this chapter provides such an expansion.}\label{Jfig4}
\end{figure}

If we know $\Ga$ and $\Gb$ (or bounds for them) then it makes sense to
use a transformation which maps the complex plane minus the slit
$[-\Gb,-\Ga]$ to the unit circle $|z|=1$. Since the transformation
\beq
t = \frac{(\Gs_1+\Ga)(1+\Gb)}{(\Gs_1+\Gb)(1+\Ga)} =
1+\frac{(\Gs_1-1)(\Gb-\Ga)}{(\Gs_1+\Gb)(1+\Ga)}
\eeq{J3}
maps the interval $[-\Gb,-\Ga]$ to $[-\infty,0]$ and takes $\Gs_1=1$ to
$t=1$, it is clear that the transformation
\beq
z = \frac{\sqrt{t}-1}{\sqrt{t}+1}\quad\text{(with $t$ given above)}
\eeq{J4}
maps the complex $\Gs_1$ plane minus the slit $[-\Gb,-\Ga]$ to the unit
disk $|z|=1$ (see \figref{Jfig4}).
So we want to get a Fast Fourier Transform (FFT) method%
\index{Fast Fourier Transform methods}
associated with
an expansion
\beq
\frac{\Gs_*}{\sqrt{[(\Gs_1+\Ga)(1+\Gb)]/[(\Gs_1+\Gb)(1+\Ga)]}} = 1+\sum_{n=1}^\infty c_n
\left(
\frac{\sqrt{\frac{(\Gs_1+\Ga)(1+\Gb)}{(\Gs_1+\Gb)(1+\Ga)}}-1}
{\sqrt{\frac{(\Gs_1+\Ga)(1+\Gb)}{(\Gs_1+\Gb)(1+\Ga)}}+1}
\right)^n.
\eeq{J5}
It looks rather formidable but really it is just a matter of substituting
\beq
t = \frac{(\Gs_1+\Ga)(1+\Gb)}{(\Gs_1+\Gb)(1+\Ga)}
\eeq{J6}
in the \cite{Eyre:1999:FNS} scheme associated with
\beq
\Gs_*/\sqrt{t} =1+ \sum_{n=1}^\infty c_n\left(\frac{\sqrt{t}-1}{\sqrt{t}+1}\right)^n.
\eeq{J7}
But of course one wants to do this substitution at the level of the
underlying Hilbert space, not just in the conductivity function.

\section{Substitution at the level of the Hilbert space}
\begin{figure}[!ht]
\centering
\includegraphics[width=0.95\textwidth]{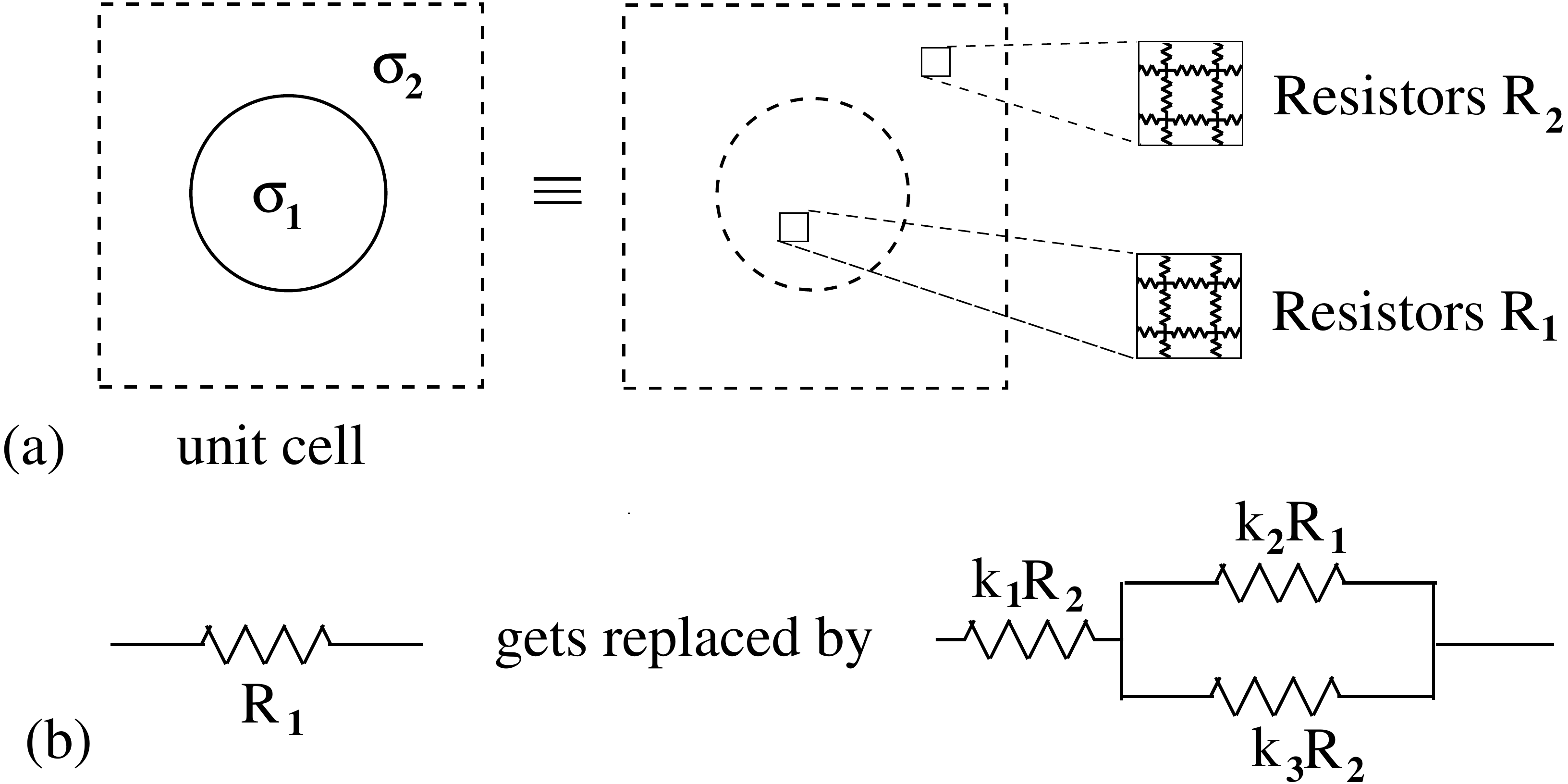}
\caption[A two-phase composite can be approximated by a discrete network of two types of resistors, and we can replace all the resistors of one type by compound resistors.]
        {A composite containing two phases, as for example a periodic array of disks with a unit cell as shown in (a), can be approximated by a discrete network
of resistors having resistors $R_2$ and $R_1$. At this discrete level we are free to replace each resistor having resistance $R_1$ by the combination of
resistors in figure (b). This transformation at the Hilbert-space level, corresponds to a fractional linear transformation in the $\Gs_1$-plane. }\label{Jfig5}
\end{figure}

\setcounter{equation}{0}
Substitution at the level of the Hilbert space
requires the substitution of subspace collections%
\index{subspace collections!substitution}
where the
subspaces do not satisfy the orthogonality property.%
\index{subspace collections!nonorthogonal}
At a discrete level
it is easy to get an idea of how this can be done. If we consider the composite
as a resistor network,%
\index{composites!approximated by resistor network}
in phase $1$ with resistors having resistance
$R_1 = \Gd/\Gs_1$, while in phase $2$ with resistors having resistance $R_2 =
\Gd/\Gs_2$, so e.g., the composite is replaced by a network (see
\figref{Jfig5}(a)).
Then at this level we could replace each resistor $R_1$ by the compound
resistor (see \figref{Jfig5}(b)) with $k_1,k_2,k_3$ real constants.
That is the resistance $R_1=\Gd/\Gs_1$ gets replaced by
\beqa
k_1R_2 + \frac{1}{1/(k_2R_1)+1/(k_3R_2)} &=&
\frac{k_1}{\Gs_2}+\frac{1}{\Gs_1/k_2+\Gs_2/k_3}\nonum
&=& \frac{k_1\Gs_1/k_2+k_1/k_3+1}{\Gs_1/k_2+1/k_3},
\eeqa{J8}
when $\Gs_2 = 1$, and so $\Gs_1$ gets replaced by
\beq
\frac{(\Gs_1/k_2+1/k_3)\Gd}{k_1\Gs_1/k_2+k_1/k_3+1},
\eeq{J9}
which is a fractional linear transformation%
\index{fractional linear transformation}
of $\Gs_1$. In other words,
in phase $1$ the number of fields is multiplied by $3$ (replacing one
resistor by $3$). Note that with real positive values of $k_1$, $k_2$ and $k_3$
the transformation \eq{J9} maps the  nonnegative real $\Gs_1$-axis onto an interval
on the positive real axis, and the negative  real $\Gs_1$-axis onto its complement, which
cannot be the desired interval $[-\beta,-\alpha]$. Essentially, instead of making
a transformation which has the effect of lengthening the branch cut as desired in the map at the top of \figref{Jfig4}, 
we have a transformation which presumably (in some suitably defined metric) shortens the branch cut.
This is why it is necessary to substitute
nonorthogonal subspace collections, rather than orthogonal ones [see also the discussion
below \eq{J24}], since, as we will see, they only act to lengthen the branch cut.

How to do this in general in the continuum case is described in section
29.1 of \cite{Milton:2002:TOC}, bottom of page 621 and page 622, for orthogonal subspace collections,
and in the second half of Section 7.8 of  \cite{Milton:2016:ETC}, without assuming orthogonality. That is a bit abstract
so let us go through it for the case in question.
\section{The original subspace collection}
Our starting point is the Hilbert space $\CH$ of
fields $\BP(\Bx)$ (real or complex $d-$dimensional vector fields), that are cell-periodic and square-integrable in the unit
cell. We denote $\BGc$ as the projection onto all fields which are zero in
phase $2$, $\CU$ as the space of constant vector fields, $\CE$ as the
space of gradients which derive from periodic potentials (i.e.,
$\BE\in\CE$ if $\nabla\times \BE = 0$ and $\langle\BE\rangle = 0$) and
$\CJ$ as the space of divergence free fields with zero average value
(i.e., $\BJ\in\CJ$ if $\nabla\cdot\BJ = 0$ and $\langle \BJ\rangle = 0$).
In our Hilbert space the inner product%
\index{inner product}
is taken to be
\beq
(\BP,\wt{\BP}) = \int_{\text{unit cell}}\overline{\BP(\Bx)}\cdot\wt{\BP}(\Bx),
\eeq{J10}
where the overline denotes complex conjugation.
With respect to this inner product the $3$ spaces $\CU,\CE$ and
$\CJ$ are orthogonal. The projection onto $\CE$ we denote by $\BGG_1$.
In Fourier space
\beq
\BGG_1\wh{\BP}(\Bk) =
\left\{\begin{array}{lr}\frac{\Bk\otimes\Bk}{|\Bk|^2}\wh{\BP}(\Bk), &
\Bk\neq 0,\\
0, & \Bk = 0.\end{array}\right.
\eeq{J11}
The field equations are solved once we have found electric fields
$\BE(\Bx)$ and current fields $\BJ(\Bx)$ such that
\beq
\BJ =
[(\Gs_1-\Gs_2)\BGc+\Gs_2]\BE,\quad\BE\in\CU\oplus\CE,\quad\BJ\in\CU\oplus\CJ.
\eeq{J12}
The effective tensor $\Gs_*$%
\index{effective!conductivity tensor}
is then defined via
\beq
\BGG_0\BJ = \Gs_*\BGG_0\BE,
\eeq{J13}
where $\BGG_0$ is the projection onto $\CU$. We assume the composite is isotropic so that $\Gs_*$ is a scalar,
although the method is easily generalized to anisotropic materials where the effective conductivity is a second order tensor.
\section[The vector subspace collection we substitute]{The vector subspace collection we substitute into the original subspace collection}
Now consider a $3-$dimensional subspace collection $\CH'$ consisting of
$3$ component vectors $\BP = [P_1,P_2,P_3]^T$ with inner product%
\index{inner product}
\beq
(\BP,\wt{\BP}) = \sum_{i=1}^3\overline{P}_i\wt{P}_i,
\eeq{J14}
where the overline denotes complex conjugation. The projection $\BGc' = p\otimes p$ projects onto the one dimensional
space of fields proportional to the unit vector $\Bp$ where $\Bp =
[p_1,p_2,p_3]^T$ and $p_1,p_2,p_3$ are given constants such that
$p_1^2+p_2^2+p_3^2=1$. The $p$'s could be complex but we \emph{do
not} mean $|p_1|^2+|p_2|^2+|p_3|^2 = 1$. Thus $\BGc'$ is a projection%
\index{projection operators!nonorthogonal}
but not an orthogonal projection
when the  $p$'s are complex, as then $\BGc' = p\otimes p$ is not Hermitian.
We take the following:
\beqa
\CU'\text{ is the space of fields proportional to } (1,0,0)^T,\nonum
\CE'\text{ is the space of fields proportional to } (0,1,0)^T,\nonum
\CJ'\text{ is the space of fields proportional to } (0,0,1)^T,\nonum
\CP_1\text{ is the space of fields proportional to }(p_1,p_2,p_3)^T,\nonum
\CP_2\text{ is the space of fields }(P_1,P_2,P_3)^T\nonum
\text{such that }p_1P_1+p_2P_2+p_3P_3 = 0.
\eeqa{J15}
The field equations become
\beq
\BJ' = [(t-\Gs_2)\BGc'+\Gs_2]\BE',\quad
\BE'\in\CU'\oplus\CE',\quad\BJ'\in\CU'\oplus\CJ',
\eeq{J16}
where the constant $t$ will be chosen so the associated ``effective
modulus'' is $\Gs_1$. That is
\beq
\BGG_0\BJ' = \Gs_1\BGG_0\BE',
\eeq{J17}
where $\BGG_0$ is the projection onto $\CU$, so that
\beq
J_1' = \Gs_1E_1'.
\eeq{J18}
Without loss of generality we can choose $E_1'=1, J_1' = \Gs_1$ so the
field equations become
\beq
\bpm
J_1'\\0\\J_3'
\epm = (t-\Gs_2)\underbrace{\bpm p_1^2 & p_1p_2 & p_1p_3\\ p_1p_2 &
p_2^2 & p_2p_3\\p_1p_3 & p_2p_3 &
p_3^2\epm}_{\BGc'}\bpm E_1'\\E_2'\\0\epm+\Gs_2\bpm E_1'\\E_2'\\0\epm.
\eeq{J19}
From the middle equation we get
\beq
(t-\Gs_2)p_1p_2E_1'+[(t-\Gs_2)p_2^2+\Gs_2]E_2' = 0,
\eeq{J20}
which with $E_1'=1$ gives
\beq
E_2' =
\frac{(\Gs_2-t)p_1p_2}{(t-\Gs_2)p_2^2+\Gs_2}.
\eeq{J21}
So we have
\beqa
\Gs_1 = J_1' &=& p_1^2(t-\Gs_2)+\Gs_2 -
\frac{(\Gs_2-t)^2p_1^2p_2^2}{(t-\Gs_2)p_2^2+\Gs_2}\nonum
&=&\Gs_2+\frac{p_1^2\Gs_2(t-\Gs_2)}{(t-\Gs_2)p_2^2+\Gs_2}\nonum
&=&\Gs_2+\frac{p_1^2\Gs_2}{p_2^2+\Gs_2/(t-\Gs_2)},
\eeqa{J22}
which with $\Gs_2=1$ is satisfied with
\beq
t = 1+\frac{\Gs_1-1}{p_1^2-p_2^2(\Gs_1-1)} =
1+\frac{(\Gs_1-1)(\Gb-\Ga)}{(\Gs_1+\Gb)(1+\Ga)},
\eeq{J23}
where
\beqa
\Ga &=& -1 -\frac{p_1^2}{p_2^2-1},\nonum
\Gb &=& -1-\frac{p_1^2}{p_2^2}.
\eeqa{J24}
Note that $-\Ga$ (respectively $-\Gb$) is obtained by substituting $t=0$
(respectively $t=\infty$) in \eq{J23}. Given real $\Gb>\Ga>0$ we need to
choose $p_1$ and $p_2$ so that these equations are satisfied. This will
necessitate complex solutions for $p_1$ and $p_2$ since otherwise $\Gb$
will be negative. Note that with $p_1$ and $p_2$ being
complex, $\BGc'$ is no longer Hermitian, even though it is a (non self-adjoint) projection, and so we have a subspace collection which is not orthogonal:%
\index{subspace collections!nonorthogonal}
$\CP_1$ is not orthogonal to $\CP_2$. Also from the field equations with $E_1' = 1$ we get
\beq
J_3' = p_3(t-\Gs_2)(p_1+p_2E_2'),
\eeq{J25}
i.e.,
\beq
J_3' = \frac{p_1p_3\Gs_2(t-\Gs_2)}{(t-\Gs_2)p_2^2+\Gs_2}.
\eeq{J26}
\section{The subspace collection after the substitution}%
\index{subspace collections!substitution}
Now consider the Hilbert space $\CH''$ consisting of all periodic fields
of the form
\beq
\BP''(\Bx) = \underbrace{\bpm
0\\\BS(\Bx)\\\BT(\Bx)\epm\BGc(\Bx)}_{\in\CP_1\otimes(\CE'\oplus\CJ')}+
\underbrace{\bpm \BQ(\Bx)\\0\\0\epm}_{\in\CH\otimes\CU'} = \bpm
\BQ(\Bx)\\\BGc(\Bx)\BS(\Bx)\\\BGc(\Bx)\BT(\Bx)\epm.
\eeq{J27}
Fields in $\CU''$ take the form
\beq
\Bu''(\Bx) = \bpm \Bu_0\\0\\0\epm\in\CU\otimes\CU'.
\eeq{J28}
Fields in $\CE''$ take the form
\beq
\BE''(\Bx) = \underbrace{\bpm0 \\ \BS(\Bx)\\
0\epm\BGc(\Bx)}_{\in\CP_1\otimes\CE'} +
\underbrace{\bpm\wt{\BE}(\Bx)\\0\\0\epm}_{\in\CE\otimes\CU'} = \bpm
\wt{\BE}(\Bx)\\\BS(\Bx)\BGc(\Bx)\\0\epm,
\eeq{J29}
where $\wt{\BE}(\Bx)\in\CE$. Fields in $\CJ''$ take the form
\beq
\BJ''(\Bx) = \underbrace{\bpm
0\\0\\\BT(\Bx)\epm\BGc(\Bx)}_{\in\CP_1\otimes\CJ'}+\underbrace{\bpm
\wt{\BJ}(\Bx)\\0\\0\epm}_{\in\CJ\otimes\CU'} = \bpm
\wt{\BJ}(\Bx)\\0\\\BT(\Bx)\BGc(\Bx)\epm,
\eeq{J30}
where $\wt{\BJ}(\Bx) \in\CJ$. The space $\CP_1'$ consists of all vectors
of the form
\beq
c\bpm p_1\\p_2\\p_3\epm,
\eeq{J31}
and $\CP_1''$ consists of all fields $\BP(\Bx)$ of the form
\beq
\bpm p_1\BC(\Bx)\\p_2\BC(\Bx)\\p_3\BC(\Bx)\epm \BGc(\Bx) \in\CP_1\otimes\CP_1'.
\eeq{J32}
Also $\CP_2'$ consists of all vectors of the form
\beq
c\bpm q_1\\q_2\\q_3\epm \quad\text{ where }~p_1q_1+p_2q_2+p_3q_3 = 0,
\eeq{J33}
and $\CP_2''$ consists of all fields $\BP(\Bx)$ of the form
\beq
\underbrace{\bpm \BQ_1(\Bx)\\\BQ_2(\Bx)\\\BQ_3(\Bx)\epm\BGc(\Bx) + (1-\BGc(\Bx))\bpm
\BR(\Bx)\\0\\0\epm}_{\in(\CP_1\otimes\CP_2'+\CP_2\otimes\CU')}\quad\text{
where }~ p_1\BQ_1(\Bx)+p_2\BQ_2(\Bx)+p_3\BQ_3(\Bx)
= 0.
\eeq{J34}
The inner product%
\index{inner product}
on $\CH''$ is defined to be
\beq
(\BP,\wt{\BP}) = \int_{\text{unit
cell}}[\ov{\BS(\Bx)}\cdot\wt{\BS}(\Bx)+\ov{\BT(\Bx)}\cdot\wt{\BT}(\Bx)]\BGc(\Bx)+\ov{\BQ(\Bx)}\cdot\wt{\BQ}(\Bx),
\eeq{J35}
where the overline denotes taking a complex conjugate. With this inner
product, the subspaces $\CU'',\CE''$ and $\CJ''$ are mutually
orthogonal. We define $\BGc''= (\Bp\otimes\Bp) \BGc$, i.e.,
\beq
\BGc'' \left\{\bpm
0\\\BS(\Bx)\\\BT(\Bx)\epm\BGc(\Bx)+\bpm\BQ(\Bx)\\0\\0\epm\right\} = \bpm
p_1^2\BI & p_1p_2\BI & p_1p_3\BI\\p_1p_2\BI & p_2^2\BI &
p_2p_3\BI\\p_1p_3\BI & p_2p_3\BI & p_3^2\BI
\epm \bpm \BQ(\Bx)\\\BS(\Bx)\\\BT(\Bx)\epm\BGc(\Bx),
\eeq{J36}
where $\BI$ is the $d\times d$ identity matrix which defines $\BGc''$
even if $p_1,p_2$ and $p_3$ are complex. Now the field equations are
\beq
\BJ'' = [(t-\Gs_2)\BGc'' +
\Gs_2\BI]\BE'',\quad\BE''\in\CU''\oplus\CE'',\quad\BJ''\in\CU''\oplus\CJ''.
\eeq{J37}
These are easy to solve given periodic solutions $\BJ(\Bx)$ and
$\BE(\Bx)$ to the equations in the Hilbert space $\CH$, i.e.,
\beq
\BJ = [(\Gs_1-\Gs_2)\BGc+\Gs_2]\BE,\quad \nabla\cdot\BJ =
0,\quad\nabla\times\BE = 0.
\eeq{J38}

We take (with $E'_1 = 1$)
\beq
\BE'' = \bpm \BE(\Bx)\\E_2'\BE(\Bx)\BGc(\Bx)\\0\epm,\quad \BJ''=\bpm
\BJ(\Bx)\\0 \\J_3'\BJ(\Bx)/\Gs_1\epm.
\eeq{J39}
Note that we have $\BE''\in\CU''\oplus\CE''$ and $\BJ''\in\CU''\oplus\CJ''$. Also, with
$\Gs_2 = 1$, we have
\beqa
((t-\Gs_2)\BGc''+\Gs_2)\BE'' &=& (\Bp\otimes
\Bp(t-\Gs_2)+\Gs_2])\bpm1\\E_2'\BE(\Bx)\BGc(\Bx)\\0\epm\nonum
&&+\bpm\BE(\Bx)(1-\BGc(\Bx))\\0\\0\epm\nonum
&=& \bpm J_1'\BE(\Bx)\\0\\J_3'\BE(\Bx)\epm\BGc(\Bx) + \bpm
\BE(\Bx)\\0\\0\epm(1-\BGc(\Bx))\nonum
&=&\bpm\BJ(\Bx)\\0\\J_3'\BE(\Bx)\epm\BGc(\Bx) +
(1-\BGc(\Bx))\bpm \BJ(\Bx)\\0\\0\epm = \BJ''.
\eeqa{J40}
Finally if $\BGG_0''$ is the projection onto $\CU''$ we have
\beqa
\BGG_0''\BE''&=&\bpm \langle\BE\rangle\\0\\0\epm\BGc(\Bx)
+\bpm\langle\BE\\0\\0\epm\rangle(1-\BGc(\Bx)), \nonum
\BGG_0''\BJ'' &=& \bpm\langle\BJ\rangle\\0\\0\epm\BGc(\Bx) +
\bpm\langle\BJ\rangle\\0\\0\epm(1-\BGc(\Bx)),
\eeqa{J41}
and since $\langle\BJ\rangle = \Gs_*\langle\BE\rangle$ we deduce that
\beq
\BGG_0''\BJ'' = \Gs_*\BGG_0''\BE''.
\eeq{J42}
That is, $\Gs_*$ is still the effective tensor. Now the idea is to apply, either
the basic Fast Fourier Transform scheme of Moulinec and Suquet%
\index{Fast Fourier Transform methods!Moulinec--Suquet scheme}
(\cite{Moulinec:1994:FNM}, \cite{Moulinec:1998:NMC})
to the Hilbert space $\CH''$ that is associated with the expansion
\beq \Gs_*/[(t+1)/2] = 1+\sum_{n=1}^\infty d_n\left(\frac{t-1}{t+1}\right)^n,
\eeq{J42a}
or the accelerated Fast Fourier Transform method of \cite{Eyre:1999:FNS}%
\index{Fast Fourier Transform methods!Eyre--Milton scheme}
as generalized in
section 14.9 of \cite{Milton:2002:TOC} to the
Hilbert space $\CH''$ that is associated with the expansion
\beq
\Gs_*/\sqrt{t} = \sum_n c_n\left(\frac{\sqrt{t}-1}{\sqrt{t}+1}\right)^n.
\eeq{J43}

The operator $\BGc''$ is easily evaluated in real space. The operator
$\BGG_1''$ which projects onto $\CE''$ is easily evaluated in Fourier
space since
\beq
\BGG_1''\bpm \BQ(\Bx)\\\BGc(\Bx)\BS(\Bx)\\\BGc(\Bx)\BT(\Bx)\epm = \bpm
\BGG_1\BQ(\Bx)\\\BGc(\Bx)\BS(\Bx)\\0\epm,
\eeq{J44}
where in Fourier space
\beq
\BGG_1\wh{\BQ}(\Bk) =
\left\{\begin{array}{lr}\frac{\Bk\otimes\Bk\wh{\BQ}(\Bk)}{|\Bk|^2},
&\Bk\neq 0,\\0,& \Bk=0.\end{array}\right.
\eeq{J45}
Hence the accelerated Fast Fourier Transform method of Eyre and Milton
can be directly applied in the Hilbert space $\CH''$.

\section{Proof of acceleration}
\setcounter{equation}{0}
Let us suppose $\Gs_1\in\RR$ is fixed such
that $\Gs_1>1$. The rate of convergence will be determined by the
magnitude of
\beq
|z| = \left|\frac{\sqrt{t}-1}{\sqrt{t}+1}\right|.
\eeq{J46}
Since $t>1$, when $\Gs_1>1$, the convergence will be quicker the smaller
$t$ is. Now as $\Ga$ increases,
\beq
t = 1+\frac{(\Gs_1-1)(\Gb-\Ga)}{(\Gs_1+\Gb)(1+\Ga)}
\eeq{J47}
decreases monotonically from the value
$t=1+[(\Gs_1-1)\Gb]/[\Gs_1+\Gb]$ at $\Ga=0$, to the value $1$ at
$\Ga=\Gb$. So the larger the value of $\Ga$, the faster the rate of
convergence. Similarly as $\Gb$ decreases, $t$ decreases monotonically
from the value $t=1+(\Gs_1-1)/(1+\Ga)$ at $\Gb=\infty$ to the value
of $1$ at $\Gb=\Ga$. So the smaller the value of $\Gb$, the faster the
rate of convergence. More generally when $\Gs_1$ is complex, to see the
rate of convergence one should plot the contours $|z|=c$ in the complex
$\Gs_1-$plane. It might be instructive to do this for particular values
of $\Ga$ and $\Gb$, and compare the contours with $\Ga=0,\Gb=\infty$.
\section{The numerical example of Moulinec and Suquet}
\labsect{Jsuqmol}
\setcounter{equation}{0}
This numerical example is due to Herv{\'e} Moulinec and Pierre Suquet (to be published) and compares in a model example the speeds of convergence of
results in the Hilbert spaces $\CH$ and $\CH''$ for the Fast Fourier Transform scheme first proposed by Moulinec and Suquet (\cite{Moulinec:1994:FNM}, \cite{Moulinec:1998:NMC}). Then the speeds of convergence of results are compared in the Hilbert spaces $\CH$ and $\CH''$
for the accelerated scheme of \cite{Eyre:1999:FNS} [see also the generalization in section 14.9 of \cite{Milton:2002:TOC}]. In both cases
the Fast Fourier Transform schemes converge substantially faster in the Hilbert space $\CH''$.

It is to be emphasized that while these numerical results are only for the effective tensor, the real interest in the new algorithm
is for obtaining results for the fields inside the body: accelerated rates of convergence for the effective conductivity alone could easily be obtained by
using Pad\'e approximants%
\index{Pad\'e approximants}
(\cite{Baker:1981:PAB})
or (almost equivalently)  the associated bounds which use series expansion coefficients up to a given order (\cite{Milton:1981:BTO}; \cite{McPhedran:1981:BET};
\cite{Milton:1982:CTM}: see also Chapter 27 in \cite{Milton:2002:TOC} and
Chapter 10 in \cite{Milton:2016:ETC}, and references therein). Using the new method I expect there
will be a similar acceleration of the convergence rates for the fields. Indeed, bounds on the norm of the difference between the actual field and the field obtained
by truncating the series expansion show these improved convergence rates. However this does not guarantee pointwise convergence, and in any case it needs to be
numerically explored. Pad\'e approximants methods could also be used for the fields, but this approach has the disadvantage that one needs to simultaneously store a lot of
information: not just the fields, but also the fields that appear up to a given order in the perturbation expansion for a nearly homogeneous medium.

\begin{figure}[!ht]
\centering
\includegraphics[width=0.8\textwidth]{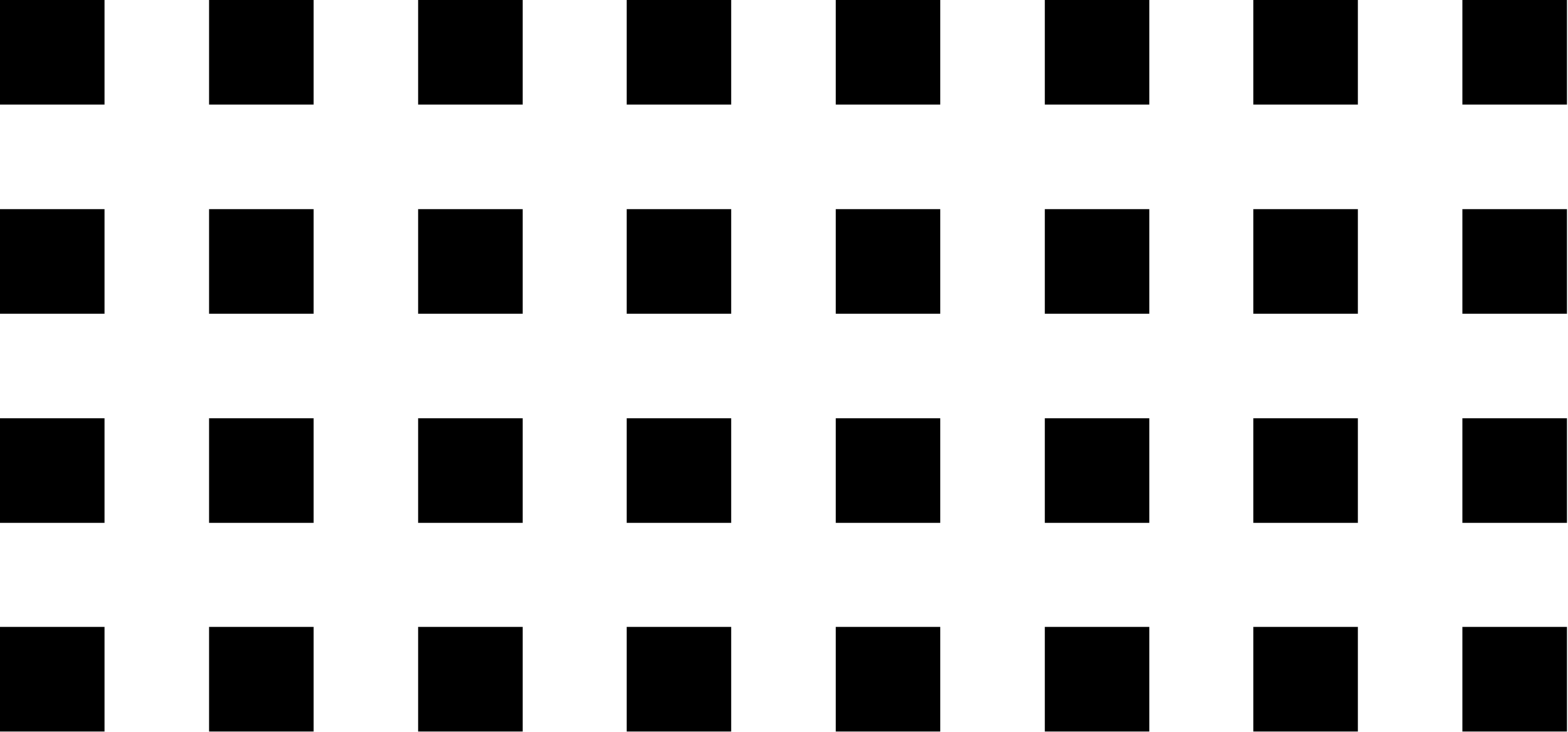}
\caption[A periodic array of squares at a volume fraction of 25\% provides a benchmark example for testing the theory.]{A periodic array of squares at a volume fraction of 25\%, as illustrated here, provides a benchmark for
testing the theory as, thanks to \protect\cite{Obnosov:1999:PHS}, there is an exact formula for its effective conductivity, given in \eq{J47a}.}\label{Jfig5a}
\end{figure}

The model example is a regular array
of squares%
\index{array of squares}
of conductivity $\Gs_1$ occupying a volume fraction of 25\% in a matrix of conductivity $\Gs_2=1$,
as illustrated in \figref{Jfig5a}.

\begin{figure}[!ht]
\centering
\includegraphics[width=0.9\textwidth]{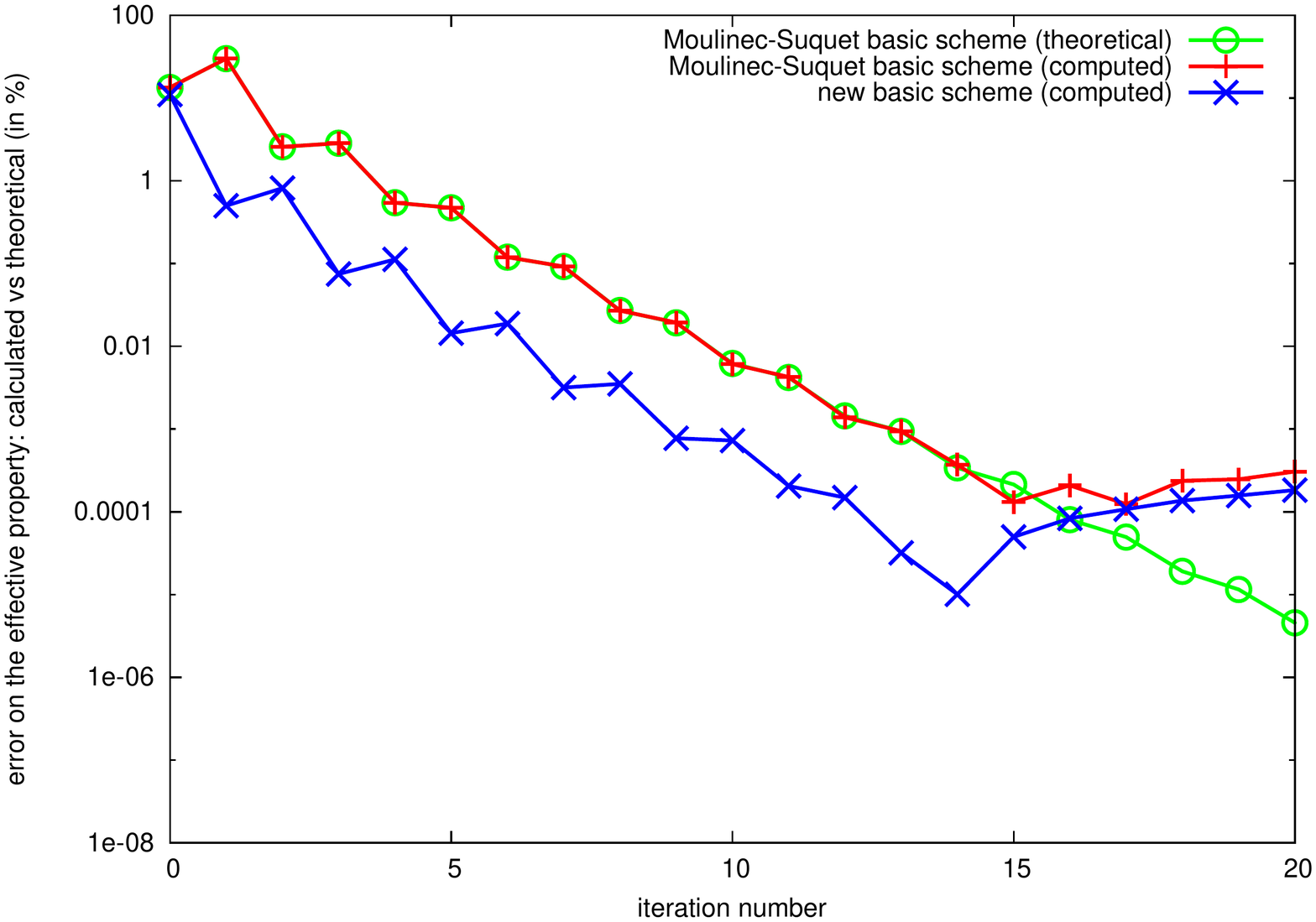}
\caption[The results of Moulinec and Suquet showing the acceleration using the new method in the rate of convergence gained over their basic scheme.]{The results of Moulinec and Suquet (to be published)
comparing the convergence of numerical
results to the formula \protect\eq{J47a} for the effective conductivity $\Gs_*$
when $\Gs_1=0$ and $\Gs_2=1$, for their
basic scheme (Moulinec and Suquet \protect\cite{Moulinec:1994:FNM}, \protect\cite{Moulinec:1998:NMC}),
and for the basic scheme applied in the Hilbert space $\CH''$. Printed with permission
of Herv{\'e} Moulinec and Pierre Suquet.}\label{Jfig5b}
\end{figure}

\begin{figure}[!ht]
\centering
\includegraphics[width=0.9\textwidth]{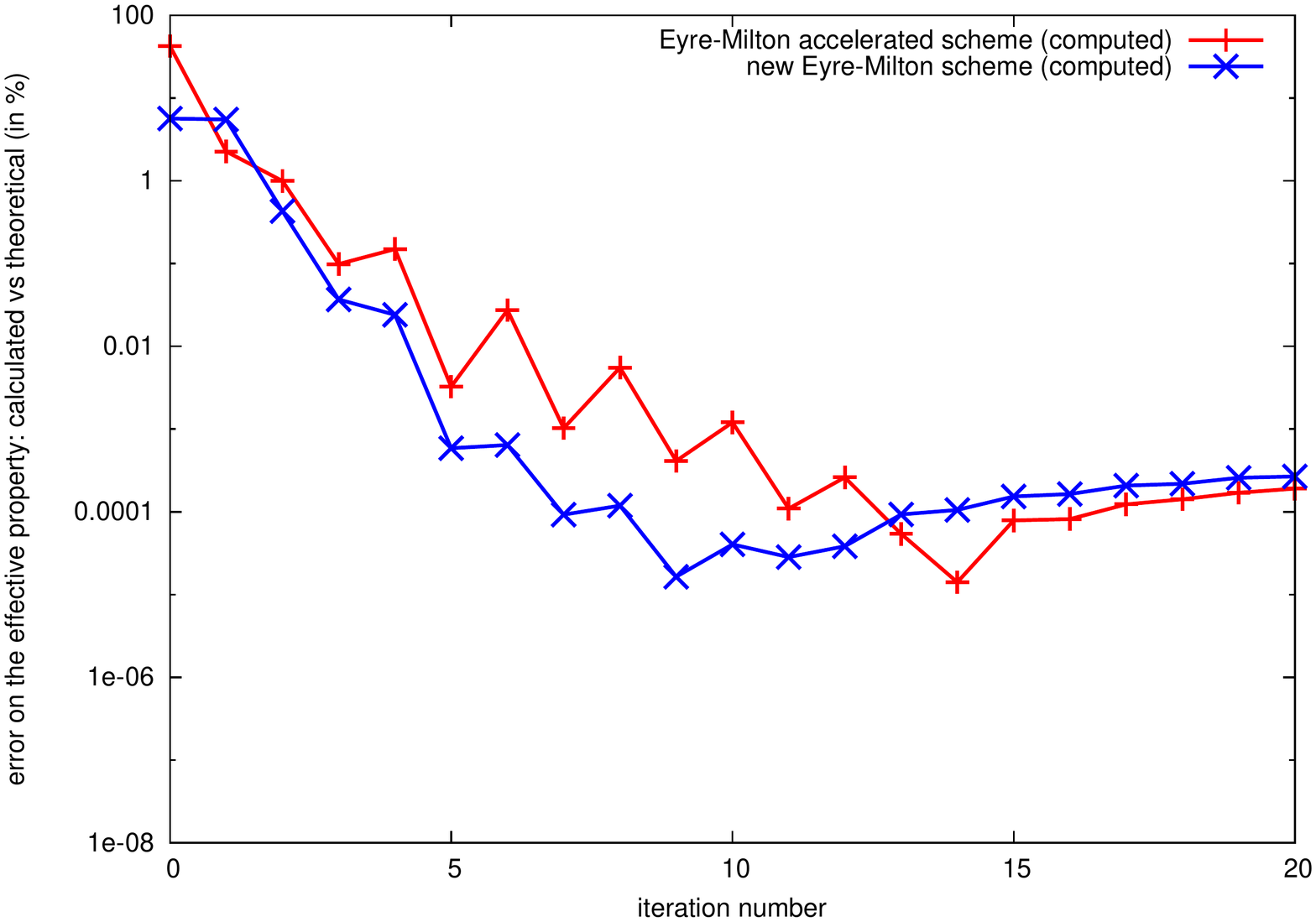}
\caption[The results of Moulinec and Suquet showing the acceleration using the new method in the rate of convergence gained over the scheme of Eyre and Milton.]
        {The results of Moulinec and Suquet (to be published) comparing the convergence of numerical
results to the formula \protect\eq{J47a} for the effective conductivity $\Gs_*$, when $\Gs_1=0$ and $\Gs_2=1$
for the accelerated scheme of \protect\cite{Eyre:1999:FNS} (see also the generalization in
section 14.9 of \protect\cite{Milton:2002:TOC}), and for the accelerated scheme applied in the Hilbert space $\CH''$. Printed with permission
of Herv{\'e} Moulinec and Pierre Suquet}\label{Jfig5c}
\end{figure}

An exact formula for the effective conductivity of this array (and for the interior fields)
was discovered by \cite{Obnosov:1999:PHS},
\beq \Gs_*=\sqrt{(1+3\Gs_1)/(3+\Gs_1)},
\eeq{J47a}
and clearly has a branch cut between $\Ga=1/3$ and $\Gb=3$. [Interestingly, for the four-phase checkerboard%
\index{checkerboard!four--phase}
\cite{Mortola:1985:TDH} conjectured a formula, that was later independently proved by
\cite{Craster:2001:FPC} and \cite{Milton:2001:PCC}.] Taking a wider estimate for the branch
cut (with $\Ga=1/4$ and $\Gb=4$), Moulinec and Suquet used the algorithm described in this chapter
and found that the acceleration provided by the new method was generally substantially improved
as shown in their \figref{Jfig5b} and \figref{Jfig5c}.

\section{Estimating the parameters $\Ga$ and $\Gb$}
\begin{figure}[!ht]
\centering
\includegraphics[width=0.75\textwidth]{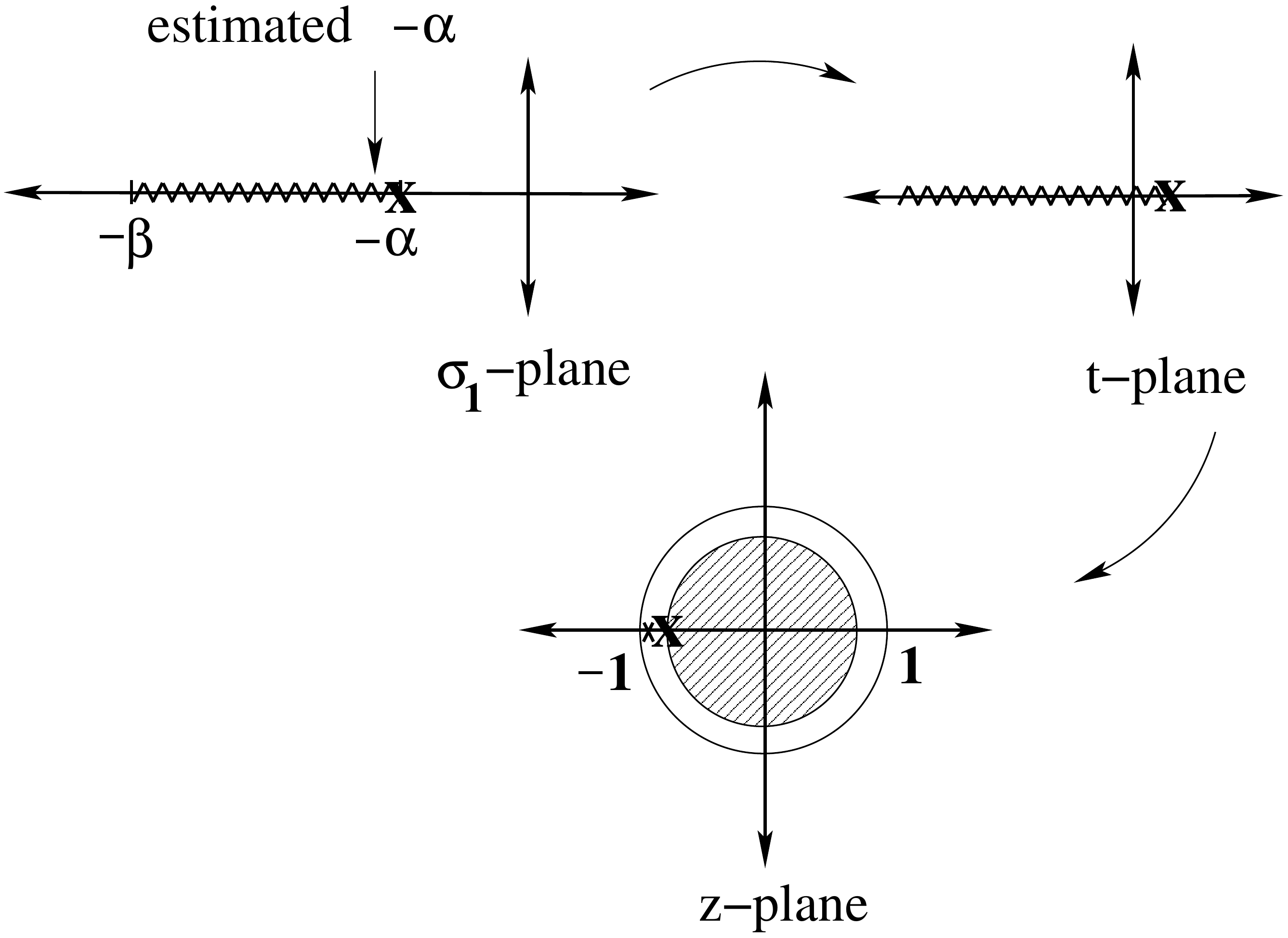}
\caption[The new accelerated scheme should still work even when we make a small error in the estimation of $\Ga$ (and/or $\Gb$).]{The new accelerated scheme should still work although not as effectively,
even when we make a small error in the estimation of $\Ga$ (and/or $\Gb$).
With this error a singularity at $\BX$, gets mapped to a point slightly greater than zero in the $t$-plane, and then gets mapped to a point
slightly inside the unit disk in the $z$-plane. There is a consequent reduction in the radius of convergence, and an associated decrease in the rate of convergence,
for the series in powers of $z$.}\label{Jfig6}
\end{figure}

\setcounter{equation}{0}
\cite{Bruno:1991:ECS} has derived rigorous lower bounds%
\index{bounds!Bruno}
on $\Ga$ and
upper bounds on $\Gb$, for suspensions of separated spheres in a medium.
Actually we need only estimates on $\Ga$ and $\Gb$. For
instance if an estimate on $\Ga$ is slightly too large, then the
transformations will look like those depicted in \figref{Jfig6}. The series will still
converge but only for $|z|<c$, where $c$ will
be close to $1$ (but less than $1$).
However, the support of the measure%
\index{measure!support}
will be highly dependent on small details of the
microstructure: a tiny sharp cusp on the surface of an otherwise smooth inclusion will in general dramatically
change the support of the measure: see \cite{Hetherington:1992:CSC}, page 378 of
\cite{Milton:2002:TOC}, and equation (15) of \cite{Helsing:2011:SSR}. As remarked on page 378 of
\cite{Milton:2002:TOC}, this is related to the fact that these sharp corners%
\index{corners as energy sinks}
behave as sinks of energy in the mathematically
equivalent dielectric problem with the conductivities $\Gs_1$ and $\Gs_2$ being replaced by electrical permittivities $\Gve_1$ and $\Gve_2$,
where $\Gve_2$ is real and positive, while $\Gve_1$ is almost real and negative with a tiny positive imaginary part.
Related observations and analysis include those of \cite{Qiu:2008:PLS}, \cite{Pikonen:2010:CFS}, \cite{Estakhri:2013:PUB},
and \cite{BonnetBenDhia:2013:RCN}. Also there is a close connection with the essential spectrum of the Neumann--Poincar{\'e} operator
on planar domains with corners (\cite{Helsing:2016:NSM}; \cite{Perfekt:2016:ESN}).

\section{Bounds on the support of the measure using $Q^*_C$-convex functions}
\labsect{spectbds}
\setcounter{equation}{0}
For many other equations in periodic composite materials robust things can be said about the measure which restrict
its support. We assume we have
a Hilbert space $\CH$ of square integrable periodic fields $\BP(\Bx)$ with an inner product
\beq (\BP',\BP)=\lang\overline{\BP}'\cdot\BP\rang, \eeq{J47aa}
where the angular brackets denote an average over the unit cell of periodicity.
We also assume $\CH$ has a decomposition
into three orthogonal subspaces
\beq \CH=\CU\oplus\CE\oplus\CJ. \eeq{J47b}
We are interested in the equations
\beq \BJ_0+\BJ(\Bx)=\BL(\Bx)(\BE_0+\BE(\Bx)), \eeq{J47c}
where
\beq \BJ_0,\BE_0\in\CU,\quad\BE\in\CE,\quad\BJ\in\CJ, \eeq{J47d}
and where the local tensor field $\BL(\Bx)$ takes the form
\beq \BL(\Bx)=\BQ(\BR(\Bx))\left[\sum_{i=1}^n\Gc_i(\Bx)\BL_i\right][\BQ(\BR(\Bx))]^T, \eeq{J47e}
in which $\BQ(\BR)$ is the orthogonal matrix (satisfying $\BQ\BQ^T=\BI$) associated with a rotation $\BR$ acting on elements in the tensor space, while $\BR(\Bx)$ is a field of rotation matrices giving the local orientation of each phase,
and $\Gc_i(\Bx)$ is the indicator function%
\index{indicator function}
that is $1$ in phase $i$ and zero elsewhere.

First, following section 14.8 of \cite{Milton:2002:TOC} let us show that there are series expansions%
\index{series expansions}
for the effective tensor
that converge even when the local tensor $\BL(\Bx)$ is not everywhere positive definite. To obtain the series
expansion we introduce a constant reference tensor%
\index{reference tensor}
$\BL_0$ and an associated operator $\BGG$:%
\index{Gamma operator}
we say
\beq \BE'=\BGG\BP\quad\text{if and only if  }\BE'\in\CE\text{ and  }\BP-\BL_0\BE'\in\CU\oplus\CJ. \eeq{J47f}
Explicitly, $\BGG$ is given by the operator
\beq \BGG=\BGG_1(\BGG_1\BL_0\BGG_1)^{-1}\BGG_1, \eeq{J47g}
where the inverse is to be taken on the subspace $\CE$. Introducing the polarization field%
\index{polarization field}
\beq \BP=(\BL-\BL_0)(\BE_0+\BE)=\BJ_0+\BJ-\BL_0(\BE_0+\BE), \eeq{J47h}
we see it satisfies $\BGG\BP=-\BE$, and hence we have the identity
\beq [\BI+\BGG(\BL-\BL_0)](\BE_0+\BE)=\BE_0+\BE-\BE=\BE_0, \eeq{J47i}
which gives
\beq \BJ_0+\BJ=\BL(\BE_0+\BE)=\BL[\BI+\BGG(\BL-\BL_0)]^{-1}\BE_0. \eeq{J47j}
Averaging both sides yields a formula for the effective tensor, and associated series expansion:%
\index{effective!tensor!formula}%
\index{series expansions}
\beq \BL_*=\lang\BL[\BI+\BGG(\BL-\BL_0)]^{-1}\rang=\sum_{j=0}^\infty\lang\BL[\BGG(\BL_0-\BL)]^j\rang,
\eeq{J47k}
where the angular brackets denote a volume average.

A sufficient condition for convergence
is that the operator $\BGG(\BL_0-\BL)$ has norm less than $1$. We take a reference tensor
of the form
\beq \BL_0=\Gs_0\BI+\BL_0', \eeq{J48}
where the tensor $\BL'_0$ is assumed to be bounded and self-adjoint, but need not be positive definite. We are interested in seeing
whether the expansion for the effective tensor converges when $\Gs_0$ is very large. Expanding the operator $\BGG$ in powers
of $1/\Gs_0$ gives
\beq \BGG=\BGG_1(\BGG_1\BL_0\BGG_1)^{-1}\BGG_1=\BGG_1/\Gs_0-\BGG_1\BL_0'\BGG_1/\Gs_0^2+\BR_1, \eeq{J49}
with a remainder term
\beq \BR_1=\sum_{j=2}^\infty\BGG_1(-1)^j(\BL_0'\BGG_1)^j/\Gs_0^{j+1}. \eeq{J50}
When $\Gs_0>\|\BL_0'\|$ this remainder term has norm satisfying the bound
\beq \|\BR_1\|\leq \sum_{j=2}^\infty\|\BL_0'\|^j/\Gs_0^{j+1}=\frac{\|\BL_0'/\Gs_0\|^2}{\Gs_0-\|\BL_0'\|}.
\eeq{J51}
This gives us a bound on the norm%
\index{bounds!norm}
of  $\BGG(\BL_0-\BL)$:
\beqa \|\BGG(\BL_0-\BL) \| & = & \|(\BGG_1/\Gs_0-\BGG_1\BL_0'\BGG_1/\Gs_0^2+\BR_1)(\Gs_0\BI+\BL_0'-\BL)\| \nonum
& = & \|\BGG_1[\BI+(\BL_0'-\BL-\BL_0'\BGG_1)/\Gs_0]+\BR_2\| \nonum
& \leq & \|\BI+(\BL_0'-\BL-\BGG_1\BL_0'\BGG_1)/\Gs_0\|+\|\BR_2\|,
\eeqa{J52}
where
\beq \BR_2\equiv \BR_1(\Gs_0\BI+\BL_0'-\BL)-\BGG_1\BL_0'\BGG_1(\BL_0'-\BL)/\Gs_0^2.
\eeq{J53}
Now we can find a $\Gs^+>\|\BL_0'\|$ such that for $\Gs_0>\Gs_0^+$,
\beq \|(\Gs_0\BI+\BL_0'-\BL)\|<2\Gs_0,\quad \Gs_0-\|\BL_0'\|>\Gs_0/2. \eeq{J54}
Hence there exists a constant $\Gb_1>0$ (independent of $\Gs_0$) such that
\beqa \|\BR_2\|&\leq & \|\BR_1\|\|(\Gs_0\BI+\BL_0'-\BL)\|+\|\BGG_1\BL_0'\BGG_1(\BL_0'-\BL)\|/\Gs_0^2 \nonum
&\leq & 4\|\BL_0'/\Gs_0\|^2+\|\BGG_1\BL_0'\BGG_1(\BL_0'-\BL)\|/\Gs_0^2 \leq \Gb_1/\Gs_0^2.
\eeqa{J55}
Now given any bounded operator $\BA$ with adjoint $\BA^\dagger$ and field $\BP\in\CH$ with $|\BP|=1$ we have
\beqa |(\BI-\BA/\Gs_0)\BP|^2 & = & (\BP,(\BI-\BA/\Gs_0)^\dagger(\BI-\BA/\Gs_0)\BP \nonum
 & = & 1-(\BP,(\BA+\BA^\dagger)\BP)/\Gs_0+(\BP,\BA^\dagger\BA\BP)/\Gs_0^2 \nonum
& \leq & 1-(\BP,(\BA+\BA^\dagger)\BP)/\Gs_0+\|\BA^\dagger\BA\|/\Gs_0^2,
\eeqa{J56}
and from the definition of the norm of an operator, this implies
\beq \|\BI-\BA/\Gs_0)\|=1+\|\BA^\dagger\BA\|/\Gs_0^2-\min_{\substack{\BP\in\CH \\ |\BP|=1}}(\BP,(\BA+\BA^\dagger)\BP)/\Gs_0.
\eeq{J57}
Setting $\BA=\BL-\BL_0'+\BL_0'\BGG_1$ we see there exists a constant $\Gb_2>0$ such that for $\Gs_0>\Gs_0^+$,
\beq  \|\BGG(\BL_0-\BL) \|\leq 1+\Gb_2/\Gs_0^2-2\min_{\substack{\BP\in\CH \\ |\BP|=1}}(\BP,(\BL_S-\BL_0'+\BGG_1\BL_0'\BGG_1)\BP)/\Gs_0,
\eeq{J58}
where $\BL_S=(\BL+\BL^\dagger)/2$. If the operator $\BL_S-\BL_0'+\BGG_1\BL_0'\BGG_1$ is coercive%
\index{coercivity property}
in the
sense that there exists a constant $\Ga'>0$ such that
\beq \BL_S-\BL_0'+\BGG_1\BL_0'\BGG_1>\Ga'\BI, \eeq{J59}
then \eq{J58} implies
\beq \|\BGG(\BL_0-\BL) \|< 1+\Gb_2/\Gs_0^2-\Ga'/\Gs_0, \eeq{J60}
and this will surely be less than $1$ for sufficiently large $\Gs_0$. Furthermore suppose the constant tensor
$\BL_0'$ is positive semi-definite on the subspace $\CE$, i.e.,
the associated quadratic form
\beq f(\BA)=\overline{\BA}\cdot\BL_0'\BA \eeq{J60.a}
is $Q^*_C$-convex,%
\index{QC@$Q^*_C$-convexity}
meaning it satisfies
\beq \lang f(\BE) \rang\geq 0,\quad \text{for all  }\BE\in\CE, \eeq{J60.b}
where the angular brackets denote a volume average over the unit cell $C$ of periodicity.
Also suppose there exists a positive constant $\Ga'$ such that
\beq \BGG_1\BL_0'\BGG_1\geq 0,\quad \BL_S>\BL_0'+\Ga'\BI. \eeq{J61}
Then \eq{J59} will be satisfied. Now each term in the series expansion for $\BL_*$ will be a polynomial in the elements
of the matrices $\BL_1,\BL_2$,$\ldots$,$\BL_n$ and hence in the domain where the series expansion converges, the effective
tensor will be an analytic function%
\index{analyticity in component moduli}
of these matrix elements. Hence the support of the measure%
\index{bounds!support of the measure}
must lie outside the region where
$\BL_S>\BL_0'$.

Note that the definition \eq{J60.b} is slightly different to the meaning of $Q^*_C$-convexity given in \cite{Milton:2013:SIG} as here
we allow for the space $\CE$ to contain certain fields whose volume average is not zero (as happens
when we treat the Schr{\"o}dinger equation,%
\index{Schr{\"o}dinger equation}
with sources and with the energy having an imaginary part: see sections 13.6 and 13.7 of  \cite{Milton:2016:ETC}.
When the fields in $\CE$ satisfy homogeneous linear differential constraints, the condition \eq{J60.b} is easily reduced to an
algebraic condition by taking Fourier transforms of \eq{J60.b}. The projection onto the space $\CE$ is given by projections $\BGG_1(\Bk)$, that are local
in Fourier space and \eq{J60.b} holds if an only if
\beq \BGG_1(\Bk)\BL_0'\BGG_1(\Bk)\geq 0, \eeq{J61a}
for all $\Bk$, including $\Bk=0$ if $\CE$ contains fields whose volume average is not zero, i.e., $\BGG_1(0)\ne 0$.
As shown in \cite{Milton:2013:SIG} this is basically the same
procedure as followed by Murat and Tartar (\cite{Tartar:1979:CCA}; \cite{Murat:1985:CVH}; \cite{Tartar:1985:EFC})
to obtain algebraic conditions for a quadratic function to be quasiconvex.%
\index{quasiconvex!quadratic forms}
For quadratic functions, quasiconvexity and $Q^*$-convexity
are equivalent when the fields satisfy differential constraints involving derivatives of fixed order, but are not equivalent when the
differential constraints have derivatives of mixed orders, as in time-harmonic wave-equations.

\section*{Acknowledgments}
G.W. Milton wishes to especially thank Herv{\'e} Moulinec and Pierre Suquet for helping simplify
his notation for the fields in the space $\CH''$, and for doing the numerical experiments
reported in section \sect{Jsuqmol}, thus providing evidence that the method suggested here works in practice,
and not only in theory. Additionally G.W. Milton
is grateful to the National Science Foundation for support through grant DMS-1211359.


\end{document}